\documentclass[11pt]{amsart}

\newtheorem{theorem}{Theorem} [section]
\newtheorem{prop}[theorem]{Proposition}
\newtheorem{lemma}[theorem]{Lemma}
\newtheorem{cor}[theorem]{Corollary}

\numberwithin{equation}{section}

\newcommand\C{{\bf C}}
\newcommand\Chat { {\hat{\C}} } 

\renewcommand\P{{\bf P}}

\newcommand\del{\partial}

\newcommand\eps{\varepsilon}

\renewcommand\phi{\varphi}
\newcommand\iso{\simeq} 

\newcommand\Aut{\mathrm{Aut}}
\newcommand\SL{\mathrm{SL}}

\renewcommand\mod{\operatorname{mod}}  

\newcommand\diam {\mathrm{diam}}
\newcommand\Rat  {\mathrm{Rat}} 
\newcommand\Poly {\mathrm{Poly}} 
\newcommand\M {\mathrm{M}} 
\newcommand\Ratbar {\overline{\Rat}}  
\newcommand\Mbar {\overline{\M}}  

\newcommand\Polybar{\overline{\Poly}}
\newcommand\MP {\mathrm{MP}}
 
\newcommand\MPbar {\overline{\MP}}
\newcommand\MPhat {\widehat{\MP}}
\newcommand\Tbar {\overline{T}}

\newcommand{\cT}{{\mathcal T}}

\renewcommand\bold {\bigskip\noindent\bf }

\usepackage{xypic}

\begin{document}

\title{Finiteness for degenerate polynomials}

\author{Laura DeMarco}

\begin{abstract}
Let $\MP_d$ denote the space of polynomials $f:\C\to\C$ 
of degree $d\geq 2$, modulo conjugation by $\Aut(\C)$.  Using
properties of polynomial trees (as introduced in \cite{DM:trees}), we show that if 
$f_n$ is a divergent sequence of polynomials in $\MP_d$, then any subsequential 
limit of the measures of maximal entropy $m(f_n)$ will have finite support.  With similar 
techniques, we observe that the iteration maps $\{\MPbar_d\dashrightarrow \MPbar_{d^n}: n\geq 1\}$
between GIT-compactifications can be resolved simultaneously with only finitely many 
blow-ups of $\MPbar_d$.  
\end{abstract}


\maketitle

\thispagestyle{empty}

\section{Introduction}

\bigskip

The goal of this article is to present  two
consequences of the properties of polynomial trees, as studied in \cite{DM:trees}. 
Both can be described as {\em finiteness statements} for degenerating families
of polynomials.  
For each degree $d\geq 2$, let 
	$$\MP_d = \Poly_d/\Aut(\C)$$ 
denote the moduli space of polynomials $f(z) = a_dz^d + \cdots + a_1z + a_0$,
$a_d\not=0$, where the affine transformations
of $\C$ act by conjugation.  We look at sequences of polynomials whose
conjugacy classes diverge in $\MP_d$ and study their limiting dynamical behavior.
Neither of the two main theorems is directly related to trees, but the tree structure
provides a natural language with which to formulate the proofs.  

{\bold Maximal measures.}  
For  a polynomial $f:\C\to\C$
of degree $d\geq 2$, let $m(f)$ denote its measure of maximal
entropy  (see \cite{Brolin}, \cite{Lyubich:entropy}).  In \cite{D:measures}, we 
studied weak limits of the measures of maximal entropy for sequences of 
rational functions which diverge in the space of all rational functions, $\Rat_d$
(with the topology of uniform convergence on $\Chat$).  
Every subsequential limit of the measures has atoms, and in
the generic case, the limiting probability 
measure is expressible as a countably infinite
sum of delta masses.  In contrast with the rational setting, we show here:

\begin{theorem}  \label{finitemeasures}
For any sequence of polynomials $f_n$ of degree $d\geq 2$ which diverges in 
$\MP_d$, every subsequential limit $\mu$ of the maximal measures $m(f_n)$
has finite support with rational masses in the ring $\mathbb{Z}[1/d]$.  
\end{theorem}

The number of points, however, in the support of the limiting
measures $\mu$ cannot be bounded in terms of the degree (see \S\ref{measures}).

{\bold Moduli space compactification.}
Let $\M_d = \Rat_d/\Aut(\Chat)$ 
denote the moduli space of rational functions of degree $d$, where the 
M\"obius transformations act by conjugation,
and let $\Mbar_d$ be the 
GIT-compactification (over $\C$) as defined in \cite{Silverman}.  It is a normal, projective
variety, with boundary $\del \M_d$ of codimension 1 \cite[Theorem 2.1]{Silverman}, 
but the iteration map 
	$$\M_d\to \M_{d^n}$$	$$f\mapsto f^n$$
does not extend continuously to this boundary for any $d\geq 2$ and $n\geq 2$
\cite[Theorem 10.1]{D:moduli2}.  
To resolve the discontinuity, 
we defined $\hat{\M}_d$ to be the closure of $\M_d$ in the infinite product 
$\prod_n \Mbar_{d^n}$ via the embedding $f\mapsto (f, f^2, f^3, \ldots)$.
For degree $d= 2$,  it was shown that $\hat{\M}_2$
is not an analytic space:  infinitely many modifications (blow-ups) of 
$\Mbar_2\iso\P^2$ are required to resolve the indeterminacy of
iterate maps $\Mbar_2 \dashrightarrow \Mbar_{2^n}$
for all $n\geq 2$  \cite[Theorem 1.4]{D:moduli2}.  

Let $\MPhat_d$ denote the closure of the polynomial slice $\MP_d\subset\M_d$ in 
the infinite product $\prod_n \Mbar_{d^n}$ via $f\mapsto (f, f^2, f^3, \ldots)$.  
In contrast with the rational case:

\begin{theorem}  \label{finitelist}
For each $d\geq 2$, 
there exists $N = N(d)$ so that the polynomial slice $\MPhat_d$ embeds*
into the finite product 
  	$$\prod_{n=1}^{N} \Mbar_{d^n} $$
via $f\mapsto (f, f^2, \ldots, f^N)$.
\end{theorem}

\begin{cor}  \label{projective}  
There exists a projective compactification $\MPhat_d$ of $\MP_d$ 
for each $d\geq 2$ such that iteration $f\mapsto f^n$ extends c-analytically* to 
	$$\MPhat_d \to \MPhat_{d^n}$$
for all $d\geq 2$ and all $n\geq 2$. 
\end{cor}

\noindent
*With these methods, we prove only that the projective embedding is 
c-analytic (analytic away from the singularities and continuous across them).  
See e.g. \cite[Ch.~4, \S5]{Whitney:varieties}.
In  Proposition \ref{normal}, we mention a sufficient condition for ``c-analytic" to 
be replaced with ``analytic".  

{\bold Background.}  In addition to the properties of
polynomial trees, outlined in \S\ref{trees}, 
we will regularly use the following facts.  

An annulus of large modulus in $\C$ contains an essential round annulus of 
comparable modulus (see e.g. \cite[Theorem 2.1]{McMullen:CDR}).  
In particular, if $A_n$ is a sequence of annuli with $\mod(A_n)\to \infty$, then
at least one complementary component of $A_n$ has diameter 
shrinking to 0 in the spherical metric on $\Chat$.  

The moduli space $\MP_d$ is finitely covered by the affine space $\C^{d-1}$
which parameterizes the monic and centered polynomials by their coefficients
\cite[Ch.I \S1]{Branner:Hubbard:1}.  
For a polynomial $f$ of degree $d\geq 2$, let 
	$$G_f(z) = \lim_{n\to\infty} \frac{1}{d^n} \log^+|f^n(z)|$$
be its escape-rate function.  
The maximal escape rate 
	$$M(f) = \max\{ G_f(c):  f'(c)=0\} $$ 
defines a proper function $\MP_d \to [0,\infty)$
\cite[Proposition 3.6]{Branner:Hubbard:1}.  In particular, 
the connectedness locus $\{f: J(f)$ is connected$\} = \{f: M(f)=0\}$ in $\MP_d$
is compact.  

{\bold Acknowledgements.}    I would like to thank C.~McMullen and M.~Popa
for their help.  I would also like to mention Jack Milnor, for 
his articles \cite{Milnor:cubicpoly} and \cite{Milnor:quad} inspired my study of 
these moduli spaces.


\bigskip\bigskip
\section{Trees and polynomials}
\label{trees}

\bigskip
This section contains a summary of 
relevant definitions and facts from \cite{DM:trees}
about polynomials and their trees.
We assume that the polynomial $f$ has disconnected
Julia set; that is, at least one critical point of $f$ lies in the basin 
of infinity.  

The tree $\Tbar = \Tbar(f)$ 
is the quotient space under
 	$$\pi: \C\to \Tbar$$
which identifies all points within a connected component 
of a level set of $G_f$.  The polynomial $f$ induces a map 
	$$F:\Tbar\to\Tbar$$
and the escape-rate function induces a height function
	$$H: \Tbar \to [0,\infty)$$
such that $H\circ \pi = G_f$.  
The open subset $T = H^{-1}(0,\infty)$ of $\Tbar$ carries
a canonical simplicial structure, determined by the conditions:
\begin{enumerate}
\item	$F$ is a simplicial map,
\item	the vertices of $T$ consist of the grand orbits of the 
		branch points of $T$, and
\item	the height function $H$ is linear on each edge of $T$.
\end{enumerate}
The Julia set of $F$ is $J(F) = H^{-1}(0)$.  We denote
by $v_0$ the base of the tree, the highest vertex of $T$ with 
valence $\geq 3$.  

{\bold Generation.}
For any point $p\in T$, its generation is defined as
	$$N(p) = \min\{n\geq 1: H(F^n(p))>H(v_0)\}.$$
Similarly, if $e$ is an open
edge of $T$, we set $N(e) = N(p)$ for any $p\in e$.

{\bold Height metric.}
The tree carries a natural metric $d_T$ defined by 
 	$$d_T(v,v') = |H(v) - H(v')|$$
if $v$ and $v'$ are adjacent vertices, and such that $H$ is an isometry
on edges.  	In this metric, the distance between the base and the Julia set is
	$$d_T(v_0, J(F)) = M(f) := \max\{G_f(c): f'(c)=0\}$$

{\bold Degree function and measure.}
Every polynomial tree carries the data of a local degree function on edges 
	$$\deg: E(T) \to \mathbb{N}$$
where $\deg(e)$ is 
the degree with which $f$ maps the annulus $\pi^{-1}(e)$ to its image.  The degree
of an edge under an iterate $F^n$ is defined by 
  $$\deg(e, F^n) = \deg(e)\deg(F(e))\cdots\deg(F^{n-1}(e)).$$
The degree function determines an $F$-invariant probability 
measure $m_T$ on $J(F)$ such that
\begin{equation}  \label{measurevalue}
	m_T(J(e)) = \frac{ \deg(e,F^{N(e)})}{d^{N(e)}} \in \frac{1}{d^{N(e)}} \cdot \mathbb{Z}
\end{equation}
where $J(e)$ is the set of all $p\in J(F)$ whose unique path to $\infty$
passes through the edge $e$.  Observe that 
\begin{equation}  \label{measurebound}
	m_T(J(e)) \leq \left( \frac{d-1}{d} \right)^{N(e)} 
\end{equation}
for all edges $e$ 
below the base $v_0$, because $\deg(e)< d$ for these edges.  

The measure $m_T$ is the push forward
$\pi_* m(f)$ of the maximal measure of $f$.

{\bold Basepoints.}
Up to isometry, the dynamical system $(\Tbar(f),d_T, F)$
depends only on the location of the polynomial $f$ in
the moduli space $\MP_d$.
A polynomial $f$ itself picks out a scale at which to view the tree.  
To make this precise, let $\Delta = \{z : |z| \le 1\}$, and
let $p(f) \in \Tbar(f)$ be the unique
point in $\pi(\Delta) \subset \Tbar(f)$
at which the height function $H$ achieves its maximum.

{\bold Strong convergence.} 
Let $v_i \in V(T)$, $i \ge 0$, denote the adjacent sequence of 
vertices from the base $v_0$ to $\infty$, 
and let $T(k) \subset T$ denote the finite subtree
spanned by the vertices within combinatorial distance
$k$ from the base.
We say a sequence $(T_n,d_n,F_n)$ 
converges strongly if:
\begin{enumerate}
	\item
the distances $d_n(v_0,v_i)$ converge for $i=1,2,\ldots,d$;
	\item
$\lim d_n(v_0,F_n(v_0)) > 0$; and
	\item
for any $k>0$ and all $n > n(k)$, there is a simplicial isomorphism
$T_n(k) \iso T_{n+1}(k)$ respecting the dynamics.
\end{enumerate}

Suppose $(T_n,d_n,F_n)$ converges strongly.
Then there is a unique pointed simplicial complex $(T',v_0)$ with
dynamics $F' : T' \to T'$ such that $T_n(k) \iso T'(k)$ for
all $n > n(k)$, and the simplicial isomorphism respects the dynamics.
The assumptions also yield a pseudo-metric $d'$ on $T'$
as a limit of the metrics $d_n$.
Let $(T,d_T,F)$ be the metrized dynamical system obtained by
collapsing the edges of length zero to points, with the simplicial 
structure on $T$ chosen so that every vertex is in the grand orbit
of a point of valence $\geq 3$.

A local degree function on $T$ is defined as a limit of the degree
functions on the edges of $T_n$.   
For each $k>0$, pass to a subsequence so that the degree function
on the edges of the finite trees $T_n(k)$ stabilize to obtain a 
degree function on edges of  $T'$; this induces a degree function
on edges of $T$.

{\bold Pointed strong convergence.}
For a point 
$p\in\Tbar$ and any integer $k>0$, we denote by $p(k)$ the point in 
$T(k)$ closest to $p$.   We say that a 
sequence of pointed trees $(T_n, d_n, F_n, p_n)$  converges strongly if: 
\begin{enumerate}
\item	the sequence $(T_n, d_n, F_n)$ converges strongly; 
\item 	the distances $d_n(v_0, p_n)$ converge; and
\item	for any $k>0$ and all $n>n(k)$, there exists a simplicial isomorphism
		$T_n(k)\iso T_{n+1}(k)$ respecting the dynamics which takes $p_n(k)$
		to $p_{n+1}(k)$.  
\end{enumerate}
As in the space of trees without marked points, there exists a well-defined limit
$(T, d_T, F, p)$ for every strongly convergent sequence.  

{\bold Geometric topology on spaces of trees.}
Let $\cT_d$ denote the space of all polynomial trees of degree 
$d$, up to isometry preserving the dynamics.   
Let $\cT_{d,1}$ denote the set of all pointed trees $(T, d_T, F, p)$  
with one marked point $p \in \Tbar$, up 
to isometry respecting the dynamics and the marked point.
A sequence of trees or pointed trees converges in the geometric topology
to a given tree (or pointed tree) if every subsequence has a
strongly convergent subsequence with the same limit.  

A tree $(T,d_T,F)$ is normalized if $d_T(v_0, J(F)) = 1$. 

\begin{theorem}  \label{compact}  \cite[Theorem 1.3]{DM:trees}
In the geometric
topology, the set of normalized polynomial trees in $\cT_d$ is compact.  
\end{theorem}


\bigskip
\section{Convergence statements}

\bigskip
In this section, we present the key lemmas needed for the proofs of the 
two theorems.

{\bold Weights.}
Let $(T, d_T, F)$ be a polynomial tree.  Let  
	$$B(p,r) = \{x\in \Tbar: d_T(p,x)<r\}.$$ 
Suppose $p\in T$ and the ball $B(p, r)$ contains no vertices 
except possibly $p$ itself.  Then the number 
of connected components of $B(p,r)\setminus\{p\}$ coincides with the number
of components of $T\setminus \{p\}$ and of $T\setminus B(p,r)$.
 
The {\em weight} of a connected component $C$ of $\Tbar\setminus
B(p,r)$ is its measure $m_T(C)$.  

\begin{lemma} \label{weightconvergence}
Suppose that pointed trees $(T_n, d_n, F_n, p_n)$ converge in the 
geometric topology to $(T, d_T, F, p)$ with $p\in T$.  Let $B(p,2r)$
be a ball containing no vertices except possibly $p$.  
Then for all $n>>0$, the set $\Tbar_n\setminus B(p_n, r)$
has the same number of components as $\Tbar\setminus B(p,r)$
and the same set of weights.
\end{lemma}

\proof	
Pass to a strongly convergent subsequence.  
The number of components stabilizes because vertices
of valence $>2$ in $B(p_n,r)$
will collapse to $p$, and vertices of valence $>2$ 
outside $B(p_n,r)$ will be bounded away from $B(p_n,r)$.  
For all $n>>0$, 
the point $p_n$ has generation $N(p_n) \leq N(p)$.  The weights
of components of $\Tbar_n\setminus B(p_n,r)$
are determined by the degree function on edges
down to any generation greater than $N(p)$.  The degree functions
converge, so these weights converge.  
\qed

{\bold Divergent sequences of polynomials.}  
For the next three lemmas, 
we suppose that $f_n$ is a sequence of polynomials which diverges 
in moduli space $\MP_d$, so that
	$$M(f_n)\to \infty$$
as $n\to \infty$.   Let $(T_n, d_n, F_n, p_n)$ be the normalized pointed
trees, so $d_n(v_0, J(F_n))=1$ for all $n$ and $p_n$ is the 
basepoint $p(f_n)$.  Let 
	$$h_n(z) = \frac{1}{M(f_n)} G_{f_n}(z)$$
be the normalized escape-rate function of $f_n$.   Each proof involves
finding large annuli (of modulus comparable to $M(f_n)$)
near the basepoint. 

\begin{lemma}  \label{p_in_T}
Suppose $f_n$ is a sequence of polynomials which diverges in $\MP_d$
while the normalized, basepointed trees
$(T_n, d_n, F_n, p_n)$ converge to $(T, d_T, F, p)$ with $p\in T$.
For each $\eps>0$, the connected components $C_n$ of 
$\Tbar_n\setminus B(p_n,\eps)$ have
	$$\diam(\pi^{-1}(C_n))\to 0$$
in the spherical metric on $\Chat$.
If $C_n^\infty$ is the unbounded component of $\Tbar_n\setminus B(p_n,\eps)$,
then 
	$$\pi^{-1}(C_n^\infty) \to \{\infty\}$$
in the Hausdorff topology on closed sets in $\Chat$.  
\end{lemma}

\proof
For all $n>>0$, the component $C_n$ is separated from the basepoint $p_n$
by a segment of length $\eps/2$.  Thus, there exists an annulus of modulus
$\eps M(f_n)/2$ which separates $\pi^{-1}(C_n)$ from some point on the unit circle.
The set $\pi^{-1}(C_n^\infty)$ is separated by this annulus from the whole unit disk,
and thus 
	$$\pi^{-1}(C_n^\infty) \to \{\infty\}$$
as $M(f_n)\to \infty$.  For all other components, the annulus 
separates $\pi^{-1}(C_n)$ from both a point on the unit circle and $\infty$, so
	$$\diam(\pi^{-1}(C_n))\to 0$$
as $M(f_n)\to \infty$.
\qed

\begin{lemma} \label{p_in_boundary}
Suppose $f_n$ is a sequence of polynomials which diverges in $\MP_d$
while the normalized, basepointed trees
$(T_n,d_n, F_n, p_n)$ converge to $(T, d_T, F, p)$ with $p\in J(F)$.
For every $r>0$ ,  
 	$$\pi^{-1}(\Tbar_n\setminus B(p_n, r)) \to \{\infty\}$$
in the Hausdorff topology on closed sets in $\Chat$.  
\end{lemma}

\proof
Assume the sequence converges strongly.  Let $e$ be an edge of $T$ along the 
path from $p$ to $\infty$ contained in $B(p,r/3)$.  Choose $k>0$ large enough
so that $T'(k)$ contains $e$.   Then for all $n>>0$, there is 
an edge $e_n$ in $T_n$ identified with $e$ under the simplicial isomorphism
$T_n(k) \iso T'(k)$, of length $l_n(e_n) > l(e)/2$, such that $e_n$ lies
in $B(p_n, r/2)$ and $e_n$ is on the path from $p_n$ to $\infty$.

Therefore, there is an annulus of modulus $M(f_n) l_n(e_n)$ separating the 
whole unit disk from $\pi^{-1}(\Tbar_n\setminus B(p_n, r))$.  Since $M(f_n)\to
\infty$ and $l_n(e_n)\to l(e)>0$, we can conclude that the sets 
$\pi^{-1}(\Tbar_n\setminus B(p_n, r))$ converge to $\{\infty\}$ in $\Chat$. 
\qed

\begin{lemma}  \label{p_escapes}
Suppose $f_n$ is a sequence of polynomials which diverges in $\MP_d$
while the basepoints $p_n = p(f_n)$ 
in the normalized trees $(T_n, d_n, F_n)$ satisfy $d_n(v_0, p_n)\to \infty$ as $n\to \infty$.
Then for every $M\geq 0$, 
	$$\diam (h_n^{-1}([0,M])) \to 0$$
as $n\to \infty$ in the spherical metric on $\Chat$.  In particular, 
	$$\diam (K(f_n)) \to 0$$
where $K(f_n)$ is the filled Julia set of $f_n$.  
\end{lemma}

\proof
Fix $M\geq 0$ and set $M' = \max\{M,1\}$.  For all $n>>0$,
$H(p_n)> M'+1$, and therefore $h_n^{-1}((M', M'+1))$ is an annulus
of modulus $M(f_n)$ separating $h_n^{-1}([0,M])$ from both a 
point on the unit circle and $\infty$.  Consequently, 
$\diam (h_n^{-1}([0,M])) \to 0$ as $M(f_n)\to\infty$.
\qed


\bigskip\bigskip
\section{Maximal measures}
\label{measures}

\bigskip
In this section, we prove that if $f_n$ is a sequence of polynomials which 
diverges in $\MP_d$, then any limit of the maximal measures $m(f_n)$ 
has finite support.   We also give an example to show that the number of points 
in the support cannot be bounded in terms of the degree.

{\bold Proof of Theorem \ref{finitemeasures}.}
Let $(T_n, d_n, F_n, p_n)$ be the normalized pointed trees associated
to the sequence $f_n$ where $p_n$ is the basepoint $p(f_n)$.  
By passing to a subsequence, we can assume that 
	$$m(f_n) \to \mu$$
weakly, and from Theorem \ref{compact} we can 
assume that the normalized trees $(T_n, d_n, F_n, p_n)$ satisfy either
\begin{enumerate}
\item	$(T_n, d_n, F_n, p_n) \to (T,d_T, F, p)$ in the geometric topology, with $p\in T$; 
\item	$(T_n, d_n, F_n, p_n) \to (T,d_T, F, p)$ in the geometric topology, with $p\in J(F)$; or
\item	$d_n(v_0, p_n) \to \infty$,
\end{enumerate}
where $(T,d_T, F)$ is itself a normalized polynomial tree.  

Suppose we are in case (1).  
Fix $\eps>0$
so that the the ball $B(p,\eps)$ contains no vertices except 
possibly $p$.  By passing to a further subsequence, it follows from Lemmas
\ref{weightconvergence} and \ref{p_in_T} that $\mu$ has the form 
	$$\mu = \sum_C m_T(C) \delta_{z(C)}$$
for some points $z(C)\in\Chat$, 
where we sum over the connected components $C$ of $\Tbar\setminus\{p\}$.
The measure $\mu$ has finite support because the number of components
is finite.    
If $p$ has generation $N(p)$, then 
	$$\mu(\{z\}) \in \frac{1}{d^{N(p)}} \cdot \mathbb{Z}_{\geq 0} $$
for every $z\in\Chat$, from (\ref{measurevalue}).  

Suppose we are in case (2).  
From Lemma \ref{p_in_boundary}, we deduce that 
	$$m(f_n) \to \mu = \delta_\infty.$$
because $\mu(\{\infty\})\geq 1-\eps$
for any $\eps>0$ from (\ref{measurebound}).  

Suppose finally we are in case (3).
Applying Lemma \ref{p_escapes} and passing to a subsequence, we see 
that
	$$m(f_n) \to \delta_z$$
for some point $z\in\Chat$.  
\qed

{\bold Unbounded support.}
The number of points in the support of the limiting measures $\mu$ cannot 
be bounded in terms of the degree.  Consider, for example, the 
cubic polynomials 
	$$f_\eps(z) = \eps z^3 + z^2$$
as $\eps\to 0$.  
These polynomials have a fixed critical point at the origin, and
for $\eps$ small, $f_\eps$ is polynomial-like of degree 2 in a neighborhood of  
the unit disk.  In fact, $f_\eps \to z^2$ locally uniformly on $\C$ as $\eps\to 0$.  

Let $(T_\eps, d_\eps, F_\eps)$ denote the metrized tree associated to $f_\eps$.
For $\eps$ sufficiently small, 
let $v_i$ denote a sequence of consecutive vertices converging to $\pi(0)\in J(F_\eps)$.
For all $i >>0$, it is not hard to see that the valence $val(v_i)$ is given by $2 \, val(v_{i-1}) - 2$,
and thus $val(v_i)\to\infty$ as $i\to\infty$.  Furthermore, choosing representatives
of the conjugacy classes $[f_\eps]$ so that the basepoint $p_\eps$ lies at the vertex
$v_i$, it is possible to arrange so that the limiting measure has $val(v_i)$ points in 
its support.  

Note, however, that while the number of points in the support is unbounded, the total mass
remaining in $\C$ is controlled.  From inequality (\ref{measurebound}), we deduce 
that as the generation of the limiting basepoint in the tree increases, the mass 
lying in the plane tends to 0.  

Similar examples can be constructed in every degree;  for example, $f_\eps(z) = \eps z^d +
z^2$.   See also \cite[\S7]{D:measures} and compare to Corollary 7.2 there, which states
that for ``most" degenerating families of polynomials of degree $d$, the number of
points in the support of the limiting measure is bounded by $d$.


\bigskip\bigskip
\section{Algebraic limits}

\bigskip

Let $\Poly_d$ denote the space of all polynomials 
	$$f(z) = a_dz^d + a_{d-1} z^{d-1} + \cdots + a_0$$
of degree $d$.  Parametrizing by the coefficients, we find
	$$\Poly_d \iso \C^*\times\C^d.$$
Let $\Polybar_d = \P^{d+1}$ denote the compactification of $\Poly_d$
in these coordinates;  that is, each point $(a_d: a_{d-1}: \cdots : a_0 : b)
\in\P^{d+1}$ determines a pair of homogeneous polynomials, up to scale,
	$$(a_dz^d + a_{d-1}z^{d-1}w + \cdots + a_0 w^d : b w^d),$$
and the boundary of $\Poly_d$ in $\Polybar_d$
is the reducible hypersurface $\{a_d\, b = 0\}$.  We will identify a point 
$(z:w)\in\P^1$ with $z/w\in\Chat$.  

Suppose $f_n$ is a sequence in $\Poly_d$ which converges to the point
	$$g = (P(z,w)w^k : bw^d) \in \del\Poly_d$$
in $\Polybar_d$, where $P$ is chosen so that $P(1,0)\not=0$.  Then the 
graph of $f_n$ in the product $\P^1\times\P^1$ converges (in the Hausdorff
topology on closed subsets) to the zero set of the homogeneous polynomial 
	$$Q((z,w), (x,y)) = P(z,w) w^k y - b w^d x .$$
In fact, if we define holomorphic $G:\Chat\to\Chat$ by
	$$G(z:w) = \left\{ \begin{array}{ll}
			 (P(z,w):bw^{d-k}) & \mbox{ if } b\not=0  \\
			 (1:0)	& \mbox{ if } b=0
				\end{array} \right.$$ 
then $f_n\to g$ in $\Polybar_d$ if and only if:
\begin{enumerate}
\item	the zeroes $Z(f_n)$ converge (with multiplicities) 
		to the zeroes of $P$ and $\infty$, and 
\item	the polynomials $f_n$ converge locally uniformly 
		to $G$ on $\C\setminus Z(P)$.
\end{enumerate}
See \cite{D:measures} for details.  

{\bold Zeroes.}  Let $f$ be a polynomial of degree $d$ with 
disconnected Julia set, and let
$(T, d_T, F)$ be its normalized tree so that $d_T(v_0, J(F))=1$.
Let $p(f)$ be the basepoint of $f$ in $\Tbar$.   

\begin{lemma}  \label{boundedzero}
Assume $p = p(f)\in T$, and suppose $B(p,2\eps)$ contains no vertices
except possibly $p$ itself.  
For every connected component $C$ of $\Tbar\setminus B(p,\eps)$
and all $n\geq N(p)$, the set $\pi^{-1}(C)\subset\C$ contains 
exactly $m_T(C) d^n$ zeroes of $f^n$.  
\end{lemma}

\proof
For each bounded component $C$ of $\Tbar\setminus B(p,\eps)$, there is 
an edge $e$ with $N(e)=N(p)$ and $C\cap J(F) = J(e)$; thus,
	$$m_T(C) = m_T(J(e)) = \frac{\deg(e, F^{N(p)})}{d^{N(p)}} = \frac{\deg(e, F^n)}{d^n}$$
for all $n\geq N(p)$.  By construction, for each $n\geq N(p)$, the iterate
$f^n$ maps $\pi^{-1}(C)$ properly to its image with degree $\deg(e,F^n)
= d^nm_T(C)$,
and its image contains the unit disk $\Delta$.  In particular, $\pi^{-1}(C)$ 
contains exactly $d^n m_T(C)$  zeroes of $f^n$.  

By the definition of the basepoint $p = p(f)$, it follows that 
$F^n(B(p,\eps))$ is disjoint from $\pi(\Delta)$ for all 
$n\geq 1$, and therefore, all other zeroes of $f^n$ must be contained in 
$\pi^{-1}(C^\infty)$, where $C^\infty$ is the unbounded component of 
$\Tbar\setminus B(p,\eps)$.  
\qed

\bigskip
For any point $p\in\Tbar$ and each $N>0$, let $p(N)\in T$ be the closest 
point to $p$ at height $\geq 1/d^N$.  Choose $\eps>0$ so that $B(p(N),\eps)$
contains no vertices except possibly $p(N)$ itself.  
Denote by $C_N$ the unbounded
component of $\Tbar\setminus B(p(N),\eps)$.  

\begin{lemma}  \label{unboundedzero}
For $p = p(f)\in \Tbar$ and every $N>0$, the set 
$\pi^{-1}(C_N)$ contains exactly $m_T(C_N) d^N$ zeroes of $f^N$.  
\end{lemma}

\proof
For every $N>0$, the image of the origin $\pi(0)$ is contained in 
a bounded component of $\Tbar\setminus\{p(N)\}$.  
For $H(p)\geq 1$, $p=p(N)$ is the unique point at height
$H(p)$, and $m_T(C_N)=0$.  The images of $C=C_N$ under
the interates of $F$ never 
intersect $\pi(0)$, and therefore, $\pi^{-1}(C)$ contains no zeroes
of $f^N$ for any $N$. 

For $H(p)<1$, each of the connected components $C$ of 
$C_N\cap\{x\in\Tbar: H(x) \leq H(p(N))\}$ has positive measure, 
and there exists an edge $e$ such that $J(e) = C\cap J(F)$
for each $C$.
The iterate $f^N$ maps $\pi^{-1}(C)$ properly 
to its image with degree
$\deg(e, F^N) = m_T(J(e)) d^N = m_T(C)d^N$.  The image 
contains $0$ because $F^N(C)$ contains all points below the base 
$v_0$.  Therefore, $\pi^{-1}(C_N)$ contains exactly $m_T(C_N) d^N$
zeroes of $f^N$.  
\qed
\medskip

{\bold Divergent sequences of polynomials.}
As before, given $p\in\Tbar$ and an integer $N>0$, 
$C_N$ denotes the unbounded component $\Tbar\setminus B(p(N),\eps)$,
where $p(N)$ is the closest point to $p$ of height $\geq 1/d^N$ and $\eps$
is chosen small enough that $B(p(N), \eps)$ contains no vertices except 
possibly $p(N)$ itself.  

\begin{lemma}  \label{kbound}
Suppose the sequence $f_n$ diverges in $\MP_d$ while the  
normalized, basepointed trees $(T_n, d_n, F_n, p_n)$ converge to 
$(T, d_T, F, p)$.  Assume also that 
	$$f_n^N \to g_N = (P(z,w) w^k: b w^{d^N})$$
in $\Polybar_{d^N}$ for some $N$ with $P(1,0)\not=0$.  Then 
 	$$k\geq m_T(C_N) d^N.$$
\end{lemma}

\proof
We need to prove that at least $m_T(C_N) d^N$ zeroes of $f_n^N$ 
converge to $\infty$ in $\Chat$ as $n\to\infty$.   From Lemmas \ref{p_in_T} and 
\ref{p_in_boundary}, we know that for each $\eps>0$,
the unbounded components $C_n^\infty$ of $\Tbar_n\setminus B(p_n, \eps)$
satisfy $\pi^{-1}(C_n^\infty) \to \{\infty\}$ as $n\to\infty$.  For every $N$,
we have $C_{n,N}\subset C_n^\infty$, and 
Lemma \ref{unboundedzero} implies that $\pi^{-1}(C_{n,N})$
contains at least $m_{T_n}(C_{n,N}) d^N$ zeroes of $f_n^N$ (when 
$\eps$ is sufficiently small).   The pointed
trees $(T_n, d_n, F_n, p_n(N))$ converge to $(T, d_T, F, p(N))$ in the geometric
topology, and therefore $m_{T_n}(C_{n,N}) = m_T(C_N)$ for all $n>>0$ by
Lemma \ref{weightconvergence}.  
\qed

\begin{lemma}  \label{goodcase}
Suppose $f_n$ diverges in $\MP_d$ and  
normalized trees $(T_n, d_n, F_n, p_n)$ converge to 
$(T, d_T, F, p)$.  Assume also that 
	$$f_n^N \to g_N = (P(z,w) w^k: b w^{d^N})$$
in $\Polybar_{d^N}$ for some $N$ with $P(1,0)\not=0$.  
If $H(p)> 1/d^N$, then there is an assignment 
	$$C\mapsto (a_C: b_C)\in\P^1$$
of the connected components of $\Tbar\setminus\{p\}$ such that 
	$$f_n^m \to g_m = \left( \prod_C (b_C z - a_C w)^{m_T(C)d^m} : 0 \right)$$
in $\Polybar_{d^m}$ for all $m\geq N$, and $(a_C: b_C) = (1:0)$
for the unbounded component of $\Tbar\setminus\{p\}$.
\end{lemma}

\proof
The hypothesis $H(p) >1/d^N$ implies that $N(p)\leq N$, so we can apply Lemmas
\ref{boundedzero}, \ref{p_in_T}, and \ref{weightconvergence} to conclude 
that the zeroes of $f^N$ converge to 
points with multiplicities governed by the proportions $m_T(C)$.  In fact, 
this holds for $f^m$ with $m\geq N$ because $1/d^N \geq 1/d^m$.  

The polynomials $f_n^N$ are converging, uniformly away from the limiting zeroes,
to the constant $\infty$ (because the unbounded components $\pi^{-1}(C_n^\infty)$ 
are converging to $\infty$), and so we can conclude that $b=0$ in the expression 
for $g_N$.  Similiary for $f_n^m$ for all $m\geq N$.  
\qed

\begin{lemma}  \label{alg:p_escapes}
Suppose $f_n$ diverges in $\MP_d$ and normalized trees 
$(T_n, d_n, F_n, p_n)$ have $d_n(v_0, p_n)\to \infty$.  Then after 
passing to a subsequence, there exists $(a:b)\in\P^1$ such that 
	$$f_n^m \to ( (bz-aw)^{d^m}: 0 )$$
for all $m\geq 1$.
\end{lemma}

\proof
This follows from Lemmas \ref{p_escapes} and \ref{boundedzero}.  
\qed


\bigskip\bigskip
\section{Moduli space compactification}
\label{GIT}

\bigskip
Let $\MPbar_d$ denote the 
closure of the polynomial slice $\MP_d$ within the projective 
GIT-compactification $\Mbar_d$ of the moduli space of rational functions
(see \cite{Silverman}).  As in \cite{D:moduli2}, we can define 
$(\Gamma_d(n), \pi_n)$ as the blow-up of $\MPbar_d$ which resolves the 
indeterminacy of the first $n$ iterate maps $f\mapsto (f^2, f^3, \ldots, f^n)$:
$$\xymatrix{  \Gamma_d(n) \ar[d]_{\pi_n} \ar[dr] & \\ 
		\MPbar_d \ar@{-->}[r]   & \MPbar_{d^2} \times\cdots\times \MPbar_{d^n} }$$
As an analytic space, $\Gamma_d(n)$ is simply the closure of 
the graph of $f\mapsto (f^2, \ldots, f^n)$ inside
the product $\MPbar_d\times\cdots\times\MPbar_{d^n}$ and $\pi_n$ is the projection
to the first factor.  

Let $\MPhat_d$ be the inverse limit space
	$$\MPhat_d = \lim_{\longleftarrow} \Gamma_d(n)$$
where $\Gamma_d(n+1)\to \Gamma_d(n)$ is the natural projection.   
The space $\MPhat_d$ is compact and contains $\MP_d$ 
as a dense open subset.  Iteration
as a map from $\MP_d$ to $\MP_{d^n}$ extends continuously to 
	$$\MPhat_d \to \MPhat_{d^n}$$
	$$(f, f^2, \ldots) \mapsto (f^n, f^{2n}, \ldots)$$
for all degrees and all $n\geq 2$.  The extension is analytic where 
$\MPhat_d$ has the structure of an analytic space.  

We aim to show that the moduli space compactification 
$\MPhat_d$ is a projective variety 
for all $d\geq 2$.  Strictly speaking, we will only prove that
there exists $N(d) < \infty$ so that the natural projection 
\begin{equation} \label{projection}
	\Gamma_d(n) \to \Gamma_d(N(d))
\end{equation}
is an analytic homeomorphism for all $n\geq N(d)$.  In this way,
we can view $\MPhat_d$ as c-analytically embedded in the finite product 
$\MPbar_d\times \cdots\times \MPbar_{d^{N(d)}}$, which is itself projective.
Without further information on the structure of $\MPbar_d$ 
and $\Gamma_d(n)$ for every $n\geq 2$, however,
it cannot be said if the projections (\ref{projection}) are analytic 
isomorphisms for all $n\geq N(d)$.   See Proposition \ref{normal}.

{\bold GIT stability conditions.}  
Every element in $\MPbar_d$ is represented by a 
stable or semistable element in $\Polybar_d
\subset \Ratbar_d\iso \P^{2d+1}$, with respect to the conjugation action 
of $\SL_2\C$, as computed in \cite{Silverman}.  The
numerical stability criteria for points in $\Ratbar_d$ reduce to the following
for points in $\Polybar_d$  \cite[Prop 2.2]{Silverman} (see also
\cite[\S3]{D:moduli2}):  

\medskip\noindent
If the degree $d$ is even, then a point $g = (P(z,w): bw^d)\in\Polybar_d$
is stable if and only if it is semistable if and only if 
\begin{enumerate}
\item	$\deg P(z,1) > d/2$, and 
\item	if $b=0$, then the multiplicity of each zero of $P(z,1)$ is $\leq d/2$.
\end{enumerate}
If the degree $d$ is odd, then a point $g = (P(z,w): bw^d)\in\Polybar_d$
is stable (respectively, semistable) if and only if 
\begin{enumerate}
\item	$\deg P(z,1) > (d+1)/2$, ($\geq (d+1)/2$), and  
\item	if $b=0$, then the multiplicity of each zero of $P(z,1)$ is $< (d+1)/2$,
		($\leq (d+1)/2$).  
\end{enumerate}
If a point is neither stable nor semistable, it is said to be unstable.

{\bold Proof of Theorem \ref{finitelist}.}
We show that there exists an $N=N(d)$ such that
the projection (\ref{projection}) is an analytic 
homeomorphism for all $n\geq N$.  
It is analytic and surjective by construction, and so 
it suffices to prove injectivity:  {\em i.e.} 
every sequence $g = (g_1, g_2, \ldots )$ in the boundary 
	$$\del\MP_d \subset \MPhat_d \subset \prod_{n=0}^\infty \MPbar_{d^n}$$
is uniquely determined by the finite list $(g_1, g_2, \ldots, g_N)\in \MPbar_d\times 
\cdots\times \MPbar_{d^N}$.  
Consequently, the inverse limit space $\MPhat_d$ will be identified 
with $\Gamma_d(N)$ which is a subvariety of the finite product
$\MPbar_d\times\cdots\times\MPbar_{d^N}$.

We proceed in steps.  

\begin{enumerate}

\bigskip
\item
Fix $N=N(d)$ so that $$\left( \frac{d-1}{d} \right)^{N-1} < \,\frac12 \,. $$  

\bigskip
\item	
Let $f_n$ be a sequence converging to $g = (g_1, g_2, \ldots)$ 
in $\MPhat_d$.  Choose representatives in $\Poly_d$ so that 
	$$f_n^N \to g_N$$
in $\Polybar_{d^N}$ with $g_N$ semistable.  Write
	$$g_N = (P(z,w) w^k : bw^d)$$
with $P(1,0)\not=0$.

\bigskip
\item
Let $(T_n, d_n, F_n)$ be the normalized tree for $f_n$ and set $p_n = p(f_n)$
to be its basepoint.  The normalized
heights $H_n(p_n)$ remain bounded:  if there were a subsequence such 
that $H_n(p_n)\to\infty$, then Lemma \ref{alg:p_escapes} implies that 
$g_N = ((bz-aw)^{d^N}:0)$ for some $(a:b)\in\P^1$ which is an unstable 
configuration.   Therefore there is a subsequence so that 
	$$(T_n, d_n, F_n, p_n) \to (T, d_T, F, p)$$
in the geometric topology.

\bigskip
\item
We show that $H(p) > 1/d^{N-1}$.   If not, then by the choice of $N$, the 
unbounded component $C_N$ of $\Tbar\setminus\{p(N)\}$ would
have $m_T$-measure $\geq 1 - (d-1)^{N-1}/d^{N-1} > 1/2$
by (\ref{measurebound}).  
From Lemma \ref{kbound}, we have
	$$ k\geq m_T(C_N) d^N > d^N/2$$
which implies that $g_N$ is unstable.  

\bigskip
\item
We are in the setting of Lemma \ref{goodcase}, so we have 
	$$g_N = \left( \prod_C (b_C z - a_C w)^{m_T(C)d^N} : 0 \right)$$
and 
	$$f_n^m \to g_m = \left( \prod_C (b_C z - a_C w)^{m_T(C)d^m} : 0 \right)$$
for all $m\geq N$.

\bigskip\item
Suppose $d$ is even.  
The stability of $g_N$ implies that 
	$$\sum_{C\,:\, (a_C:b_C)=(1:0)} m_T(C) d^N < \frac{d^N}{2}$$
and 
	$$\sum_{C\,:\, (a_C:b_C)=(a:b)} m_T(C) d^N \leq \frac{d^N}{2}$$
for all $(a:b)\not= (1:0)$.   The same inequalities are satisfied
for every $m$ in place of $N$, so $g_m$ is stable for all $m\geq N$.

\bigskip\item
Suppose $d$ is odd.  
The semistability of $g_N$ implies that 
	$$\sum_{C\, :\, (a_C:b_C)=(1:0)} m_T(C) d^N \leq \frac{d^N-1}{2}$$
and 
	$$\sum_{C\,:\, (a_C:b_C)=(a:b)} m_T(C) d^N \leq \frac{d^N+1}{2}$$
for all $(a:b)\not= (1:0)$.  By our choice of $N$, $H(p) > 1/d^{N-1}$ implies
that $N(p) \leq N-1$.  Therefore, 
	$$m_T(C) \in \frac{1}{d^{N(p)}} \mathbb{Z} \subset \frac{1}{d^{N-1}} \mathbb{Z}$$
implies that $d$ divides $d^N m_T(C)$ for every component $C$.  But 
the largest integer divisible by $d$ and $\leq (d^N+1)/2$ is in fact 
$< (d^N+1)/2$, and therefore, 
	$$\sum_{C: (a_C:b_C)=(a:b)} m_T(C) d^N < \frac{d^N}{2}.$$
This inequality remains satisfied for all $m$ in place of $N$, and therefore
$g_m$ is stable for all $m\geq N$.  

\bigskip\item
The stability of the limit point $g_m\in\Polybar_{d^m}$ implies that $g_m$
is a representative of the $m$-th entry of $g$ for all $m\geq N$, so the $m$-th
iterates of the sequence $f_n$ converge to $g_m$ in the quotient space
$\MPbar_{d^m}$.  
The convergence is independent of the sequence we started with;
therefore, all entries of $g$ have been expressed in terms of 
$(g_1, \ldots, g_N)$.  This concludes the proof that the projection (\ref{projection})
is a homeomorphism and the proof of the theorem.

\end{enumerate}
\qed

{\bold  Proof of Corollary \ref{projective}.}
Let $N(d)$ be chosen as in Theorem \ref{finitelist}.  
Fix $n$ and choose $k\geq N(d^n)$.  
By construction, iteration $\MP_d\ni f\mapsto f^n\in\MP_{d^n}$ 
extends analytically to 
	$$\Gamma_d(kn) \to \Gamma_{d^n}(k)$$
as the projection $(f,f^2,f^3, \ldots, f^{kn})$ to $(f^n, f^{2n}, \ldots, f^{kn})$.  
Postcomposing with the analytic projection $\Gamma_{d^n}(k)\to \Gamma_{d^n}(N(d^n))$
and precomposing by the c-analytic $\Gamma_d(N(d)) \to\Gamma_d(kn)$,
we deduce that iteration extends c-analytically to 
	$$\Gamma_d(N(d))\to \Gamma_{d^n}(N(d^n))$$
for all $d\geq 2$ and all $n\geq 2$.  The graphs $\Gamma_d(N(d))$ are projective.
\qed

{\bold Normality.}
We conclude by stating a sufficient condition for the projections 
(\ref{projection}) to be isomorphisms.  

\begin{prop}  \label{normal}
If the graph $\Gamma_d(N(d))$ is normal, then 
 $$\MPhat_d\iso \Gamma_d(N(d))$$ 
is a projective variety, and iteration $f\mapsto f^n$
extends analytically to $$\MPhat_d\to \MPhat_{d^n}$$ 
for all $d\geq 2$ and all $n\geq 1$.
\end{prop}

\proof
If $\Gamma$ is normal, then any modification (blow-up) 
$\pi: X\to\Gamma$ 
which is a bijection is in fact an isomorphism.  This is a consequence
of Zariski's Main Theorem (\cite[Ch.~V, Theorem 5.2]{Hartshorne} 
applied to the inverse of $\pi$). 
Consequently, the projections 
(\ref{projection}) are isomorphisms and $\MPhat_d\iso\Gamma_d(N(d))$
for all degrees $d$. 
\qed

\bigskip\bigskip

\begin{thebibliography}{BH1}

\bibitem[BH]{Branner:Hubbard:1}
B.~Branner and J.H. Hubbard.
\newblock {The iteration of cubic polynomials. {I}. {T}he global topology of
  parameter space}.
\newblock {\em Acta Math.} {\bf 160}(1988), 143--206.

\bibitem[Br]{Brolin}
H.~Brolin.
\newblock {Invariant sets under iteration of rational functions}.
\newblock {\em Ark. Mat.} {\bf 6}(1965), 103--144.

\bibitem[De1]{D:measures}
L.~DeMarco.
\newblock {Iteration at the boundary of the space of rational maps}.
\newblock {\em Duke Math. J.} {\bf 130}(2005), 169--197.

\bibitem[De2]{D:moduli2}
L.~DeMarco.
\newblock {The moduli space of quadratic rational maps}.
\newblock {\em {\em To appear}, J. Amer. Math. Soc., {\em 2006}}.

\bibitem[DM]{DM:trees}
L.~DeMarco and C.~McMullen.
\newblock {Trees and the dynamics of polynomials}.
\newblock {\em {\em Preprint, 2006}}.
\newblock {Available at http://arxiv.org/abs/math.DS/0608759}.

\bibitem[Ha]{Hartshorne}
R.~Hartshorne.
\newblock {\em Algebraic Geometry}.
\newblock Springer-Verlag, New York, 1977.

\bibitem[Ly]{Lyubich:entropy}
M.~Lyubich.
\newblock {Entropy properties of rational endomorphisms of the {R}iemann
  sphere}.
\newblock {\em Ergodic Theory Dynamical Systems} {\bf 3}(1983), 351--385.

\bibitem[Mc]{McMullen:CDR}
C.~McMullen.
\newblock {\em Complex Dynamics and Renormalization}.
\newblock Princeton University Press, Princeton, NJ, 1994.

\bibitem[Mi1]{Milnor:cubicpoly}
J.~Milnor.
\newblock {Remarks on iterated cubic maps}.
\newblock {\em Experiment. Math.} {\bf 1}(1992), 5--24.

\bibitem[Mi2]{Milnor:quad}
J.~Milnor.
\newblock {Geometry and dynamics of quadratic rational maps}.
\newblock {\em Experiment. Math.} {\bf 2}(1993), 37--83.
\newblock With an appendix by the author and Lei Tan.

\bibitem[Si]{Silverman}
J.~H. Silverman.
\newblock {The space of rational maps on $\bf {P}\sp 1$}.
\newblock {\em Duke Math. J.} {\bf 94}(1998), 41--77.

\bibitem[Wh]{Whitney:varieties}
H.~Whitney.
\newblock {\em Complex Analytic Varieties}.
\newblock Addison-Wesley Pub. Co., Reading, MA, 1972.

\end{thebibliography}

\def\cprime{$'$}

\end{document}